\documentclass[10pt]{amsart}

\pdfoutput=1

\usepackage{tikz}
\usetikzlibrary{positioning}
\usetikzlibrary{arrows}
\usepackage{amssymb,amsmath,latexsym,graphicx}
\usepackage{charter,eucal}
\usepackage{amssymb}
\usepackage{amscd}
\usepackage[all,cmtip]{xy}
\usepackage{rotating,multicol,relsize}
\usetikzlibrary{decorations.markings}

\newenvironment{Figure}
  {\par\medskip\noindent\minipage{\linewidth}}
  {\endminipage\par\medskip}

\newtheorem{theorem}{Theorem }[section]

\newtheorem{lemma}[theorem]{Lemma}
\newtheorem{observation}[theorem]{Observation}

\newtheorem{remark}[theorem]{Remark}
\newtheorem{corollary}[theorem]{Corollary}
\newtheorem{proposition}[theorem]{Proposition}
\newtheorem{principle}[theorem]{\textsc{Principle}}

\newcommand{\bt}{\begin{theorem}}
\newcommand{\et}{\end{theorem}}
\newcommand{\bmt}{\begin{maintheorem}}
\newcommand{\emt}{\end{maintheorem}}
\newcommand{\bc}{\begin{corollary}}
\newcommand{\bl}{\begin{lemma}}
\newcommand{\ec}{\end{corollary}}
\newcommand{\el}{\end{lemma}}
\newcommand{\bo}{\begin{observation}}
\newcommand{\eo}{\end{observation}}
\newcommand{\bp}{\begin{proposition}}
\newcommand{\ep}{\end{proposition}}
\newcommand{\br}{\begin{remark}}
\newcommand{\er}{\end{remark}}
\newcommand{\bpr}{\begin{principle}}
\newcommand{\epr}{\end{principle}}

\def\Aut{\mathrm{Aut}}
\def\I{\mathop{\mathrm{I}}}

\def\PG{\mathbf{PG}}
\def\PSL{\mathbf{PSL}}
\def\GL{\mathbf{GL}}
\def\C{\mathbb{C}}

\def\eop{\hspace*{\fill}$\blacksquare$}

\def\id{\mathrm{id}}

\newcommand{\PGL}{\mathbf{PGL}}
\newcommand{\F}{\mathbb{F}}

\newcommand{\mA}{\mathcal{A}}
\newcommand{\mP}{\mathcal{P}}
\newcommand{\mL}{\mathcal{L}}

\newcommand{\bI}{\mathbf{I}}

\newcommand{\M}{\mathcal{M}}

\newcommand{\mB}{\mathcal{B}}
\newcommand{\mS}{\mathcal{S}}
\newcommand{\mD}{\mathcal{D}}

\newcommand{\mG}{\mathcal{G}}

\usepackage{amsmath,amssymb,xcolor}
\usepackage[linktoc=page,colorlinks,linkcolor={red!80!black},citecolor={red!80!black},urlcolor={blue!80!black}]{hyperref}

\title[The geometry of drums]{The geometry of drums}

\subjclass[2000]{}

\author{Koen Thas}

\thanks{}

\address{Ghent University, Department of Mathematics, Krijgslaan 281, S25, B-9000 Ghent, Belgium}

\email{koen.thas@gmail.com}

\date{}

\begin{document}

\maketitle

\begin{abstract}
We introduce the new concept of D-geometry (or ``drum geometry''), which has been recently discovered by the author in \cite{KT-DRUMS} when constructing and classifying isospectral and length equivalent drums under certain constraints. We will show that any pair of length equivalent domains, and in particular any pair of isospectral domains (which makes one unable to ``hear the shape of drums'') which is constructed by the famous Gassmann-Sunada method, naturally defines a D-geometry, and that each D-geometry gives rise to such domains.
One  goal of this letter is to show that in the present theory of isospectral and length equivalent drums, many examples are controlled by finite geometrical phenomena in a very precise sense.
\end{abstract}

%\setcounter{tocdepth}{1}
%\bigskip
%{\footnotesize
%\tableofcontents
%}

%\begin{multicols}{1}

\medskip
\section{Introduction}

\subsection{}

Ever since Mark Kac's seminal paper \cite{Kac}, a massive amount of literature has been published on ``hearing''  the shape of geometrical objects, such
as planar drums and other (higher dimensional) Riemannian manifolds. Concretely, one starts with the knowledge of an eigenspectrum (with or without 
multiplicities, depending on the context) of an unknown geometrical object with respect to a differential operator (such as the Laplacian), and one tries
to derive invariants for the object from that spectrum. A very famous example is manifested through Weyl's formula \cite{Weyl}
\begin{equation}
\lim_{x \mapsto \infty}\frac{N(x)}{x^{d/2}} = (2\pi)^{-d}\cdot \omega_d\cdot \mathrm{vol}(\Omega),
\end{equation}
where $\Omega \subset \mathbb{R}^d$ is a bounded domain, $N(x)$ is the number of Dirichlet eigenvalues $\lambda$ (with multiplicities) for which $\lambda \leq x$, and 
$\omega_d$ is the volume of a $d$-dimensional ball;  
from this formula one can conclude that the eigenspectrum (with multiplicities) of the Laplace operator of a Riemannian manifold in the Dirichlet formulation
determines the perimeter and surface of the manifold. Kac's question amounts to whether 
that spectrum also determines the {\em isometry class} of the manifold, given that the latter is a real planar drum. Recall that the Dirichlet problem for the Laplacian (for a bounded domain $\mD \subset \mathbb{R}^2$ and twice-differentiable real-valued functions $f$) is

\begin{equation}
\label{Diri}
\left\{\begin{array}{ccc}
\mathlarger{\sum_{i = 1}^2\frac{\partial^2}{\partial x_i^2}}f  + \lambda f  &=  &0\\
f_{\vert \partial \mD} &=  &0
\end{array}\right.
\end{equation}

Many variations exist presently: for instance, one has  changed the geometrical objects in an extensive way (for instance to graphs | see, e.g., \cite{KTIV}) and studied similar questions for other operators (such as the Cauchy-Riemann operator and Dirac operator). Experiments have also been done on these questions (see \cite{EvePier,Moon,SriKud}), and there are a number of papers on the subject of quantum drums \cite{TolTol}, density of drums \cite{Amore1,Amore2}, etc. | the story has gotten very big.

\subsection{}

It appeared that the higher dimensional analogue (i.e., for the $d$-dimensional Laplacian in $\mathbb{R}^d$, $d \geq 3$) of that question has a {\em negative} answer: when the dimension of the manifold is at least $3$, one can construct
nonisometric manifolds which have the same complete eigenspectrum for the Laplace operator. 

Early examples of ``flat tori'' sharing the same
eigenvalue spectrum were found in 1964 by Milnor \cite{Milnor} from
nonisometric lattices of rank $16$ in $\mathbb{R}^{16}$.  Other examples of isospectral Riemannian manifolds
were constructed later, e.g. on lens spaces \cite{Ike} or on surfaces
with constant negative curvature \cite{Vig}. In 1982, Urakawa produced the
first examples of isospectral domains in $\mathbb{R}^{n}$, $n\geq 4$ \cite{Ura}.
 Work of B\'{e}rard and Besson \cite{BerBes} is essential in that paper.
Later, in the late eighties, various other papers appeared, giving necessary conditions that any family of
billiards sharing the same spectrum should fulfill (\cite{Mel},
\cite{OsgPhiSar}, \cite{OsgPhiSar2}), and necessary conditions given as
inequalities on the eigenvalues were reviewed in \cite{Pro}.

Much more details can be found in the recent review \cite{GirTha}. 

In the original setting, which appeared to be the hardest, the first counter examples were exhibited in \cite{GWW} | see Appendix \ref{count} for more details.

\subsection{}

Call Riemannian manifolds {\em length equivalent} if they have the same set of lengths of closed geodesics, ignoring multiplicities. Note that sometimes the term ``length isospectral'' is used for the same notion. We refer to \S \ref{GS} for more details (such as related combinatorial properties).

\subsection{}

Many of the known pairs of isospectral but nonisometric manifolds arise from the {\em Gassmann-Sunada method}. (Not {\em all} do however: see e.g. the examples produced by Vign\'{e}ras in \cite{Vig} using quaternion algebras.)
This method is briefly explained in Appendix \ref{GS}; it amounts to finding {\em nonconjugate subgroups} of a finite group which are ``almost conjugate'' (AC) | see Appendix \ref{GS}. The action of the group on these subgroups by left translation explains how one unfolds a base tile to {\em make} the isospectral domains.
Then through topological considerations, one is able  to construct isospectral manifolds for which the fundamental groups act as covers of the nonconjugate subgroups.
In fact, more generally, nonisometric {\em length equivalent} manifolds can also be constructed through a variation of this method, see also Appendix \ref{GS}. To make a distinction between isospectral and length equivalent domains constructed with the Gassmann-Sunada method, we will speak of the {\em strong} Dirichlet problem in the former case, and just Dirichlet problem in the latter.
 
All the known $\mathbb{R}$-planar examples | counter examples to Kac's question | are constructed through a method called {\em transplantation}. This was long thought to be a special case of Gassmann-Sunada theory, in which one has extra information | a ``transplantation matrix'' | at hand, which explains {\em how} the tilings of the isospectral domains in the plane are connected to each other. (This is how the famous Gordon-Webb-Wolpert example is presented in Appendix \ref{count}.)
Transplantation also exists in higher dimensions.

In our recent paper \cite{KT-DRUMS}, we observed  that the transplantation method and the Gassmann-Sunada method are actually the {\em same} methods in different guises, in the following sense.

Let $(G,X)$ and $(G,Y)$ be two faithful transitive permutation groups of the same degree; one can identify 
both actions similarly with left translation actions $(G,G/U)$ and $(G,G/V)$ with $U, V$ subgroups of $G$, and $[G : U] = [G : V]$. Now $(G,X)$ and 
$(G,Y)$ have the same permutation character if and only if AC is satisfied for $(G,U,V)$. 
Then let $k \in \{ \mathbb{Q},\mathbb{R},\C \}$, and let $V := k[G/H]$ be the $k$-vector space defined by the formal sums over $k$ with base the elements of the left coset space $G/H$, $H$ any subgroup of $G$. Then $G$ acts naturally on $V$ by left translation. If $B$ is the (standard) base 
of $V$ corresponding to $G/H$, then $G$ naturally defines a subgroup $G_H$ of $\GL_{\ell}(V)$, with $\ell = [G : H]$ (by its action on $B$ and linear extension). 
In this way, there arises a faithful $k$-linear representation
\begin{equation}
\rho_H:\ \ G \ \ \hookrightarrow\ \ \GL_{\ell}(V).
\end{equation}

\bt[Transplantability \cite{KT-DRUMS}]
Two faithful permutation representations $(G,G/U)$ and $(G,G/V)$ of the same degree (with $U$ and $V$ subgroups of the group $G$) have the same permutation character if and only if AC is satisfied for $(G,U,V)$ if and only if the faithful $k$-linear representations $\rho_U$ and $\rho_V$ are $k$-linearly equivalent  if and only if there exists a $T \in \GL_{\ell}(V)$ such that for all $g \in G$, 
$\rho_U(g)\circ T = T \circ \rho_V(g)$.
\et

(Here, $k[G/U]$ and $k[G/V]$ are naturally identified.)

Buser et al. produced in \cite{Conw}  all the known planar examples of isospectral drums as domains unfolded from tiles as above, and observed that to each such domain, one can naturally associate a number of reflections of $\mathbb{R}^2$ which then generate a group. In any of the cases of \cite{Conw}, one ends up with projective general linear groups $\PGL_n(q)$ for certain low values of  $n$ and $q$. This is how all the planar (and other) examples are connected to the Gassmann-Sunada method | see also the examples in \S \ref{geomdrum}.

Using the work of Buser et al. \cite{Conw}, in \cite{Giraud}, Giraud interpreted the planar examples in {\em projective spaces}, and observed that the groups occurring in the Gassmann-Sunada triples could be seen as stabilizers of subspaces of projective spaces. In \cite{KTI,KTII,KTIII,KTIV} the author built further on the work of Giraud to settle a number of questions raised in \cite{Giraud}.

\subsection{}

Although a very large theory is at hand for the (strong) Dirichlet problem for the Laplacian, no general classification theory/method of solutions is available. Even when assuming the natural but still very general assumption of transplantability/Gassmann-Sunada,  nothing was known till recently. Let us call isospectral or length equivalent, nonisometric manifolds with respect to the Laplacian ``Gassmann-Sunada domains (or drums)'' if they can be constructed through Gassmann-Sunada triples. Again, we make a distinction between Gassmann-Sunada triples and {\em strong} Gassmann-Sunada triples in the obvious way.
The aforementioned cavity in the theory is the main concern of the author's recent paper \cite{KT-DRUMS}, in which a classification is proposed for {\em irreducible} Gassmann-Sunada domains, which could be seen as building bricks of general Gassmann-Sunada domains. We describe these ideas in some more detail in the next section. 

One of the novelties of \cite{KT-DRUMS} and the theory further described in this letter, is an intriguing connection which has surfaced between the world of irreducible drums and  that of {\em finite simple groups}, the atoms of finite group theory. The bridge between both worlds is the ``geometry of drums'' of the present letter.

For, while writing up \cite{KT-DRUMS}, it has occurred to the author that underneath all examples of Gassmann-Sunada domains for (\ref{Diri}), lies a synthetic geometry consisting
of points and lines, which satisfies a simple combinatorial property. Moreover, each such abstract geometry gives rise to  Gassmann-Sunada domains by reversing the process.  The presence of these geometries gives a satisfying explanation of the fact {\em why} the projective general linear groups play a central role in the analysis of the known planar examples, as explained in the previous subsection. (Also, they play an implicit but important role in the proofs of \cite{KT-DRUMS}.)

On the other hand, they appear to be a powerful tool for construction of domains which are isospectral or length equivalent with respect to (\ref{Diri}).  

In \S \ref{geomdrum}, we will describe the precise connection between Gassmann-Sunada drums and their geometries. We will also situate the 
projective linear groups $\PGL_m(k)$ in this context.

\subsection*{Acknowledgment}

The author wants to thank Dipendra Prasad for various highly interesting communications. 
\

\setcounter{tocdepth}{1}
{\footnotesize
\tableofcontents
}

\medskip
\section{Classifying isospectral drums}
\label{class}

Classifying isospectral drums under general constraints is obviously very hard, and almost no results on this matter are known.

In the recent paper \cite{KT-DRUMS}, a classification of {\em length equivalent} manifolds is carried through. Length equivalent manifolds are more general
than isospectral manifolds, as they correspond to manifolds with the same eigenspectrum for (\ref{Diri}) but not necessarily with the same multiplicities, 
the advantage being that inductive methods can be applied for this larger class of drum-pairs. And of course, the isospectral pairs are automatically a part of our classification.

In \cite{KT-DRUMS} we first described a number of very general classes of length equivalent manifolds (cf. the next section), in each of the constructions starting from a  given example that arises itself from the Gassmann-Sunada method. 
The constructions include the examples arising from the transplantation technique (and thus in particular the planar examples). 
To that end, we introduced four properties | called FF, MAX, PAIR and INV | inspired by natural (physical) properties (which also rule out trivial constructions), that are satisfied for each of the known  planar examples. (We refer to \cite{KT-DRUMS} for additional motivations and details about these properties.) 
Start with a Gassmann-Sunada triple $(G,U,V)$.

\medskip
$\left\{
\begin{tabular}{p{0.9\textwidth}}
\begin{itemize}
\item[{\bf FF}]
is a degeneracy property in terms of normal subgroups. If $N$ is a nontrivial normal subgroup of $G$ inside $U \cap V$, then $(G/N,U/N,V/N)$ is also 
a Gassmann-Sunada triple defining the same domains. All the necessary information is passed through to $(G/N,U/N,V/N)$.
So we want to exclude this situation, and assume that the algebraic data $(G,U,V)$ is in this sense ``minimal' with respect to the constructed domains.
\item[{\bf MAX}]
expresses the fact that the constructed domains are {\em irreducible}, by assuming that both $U$ and $V$ are {\em maximal subgroups} of $G$. If for instance $U$ is not maximal, there is a subgroup $W$ of $G$ such that $U \lneqq W \lneqq G$, and then it could happen that $V$ is also properly contained in a maximal subgroup $W'$ such that $(G,W,W')$ also is Gassmann-Sunada. If this is the case, the constructed domains from $(G,U,V)$ each consist of smaller tiles coming from $(G,W,W')$. 
The study of irreducible examples, which, in the above sense, could be seen as ``atomic examples,'' is hence the first target. 
\item[{\bf PAIR}]
expresses the fact that domains are built from the same number of tiles | this can be done in various ways, one of which being the property $\vert U \vert = \vert V \vert$. Strong Gassmann-Sunada triples automatically satisfy PAIR; by asking PAIR, we wish to study general Gassmann-Sunada triples that reflect the same behavior.
\item[{\bf INV}]
expresses the fact that the domains are guaranteed to be {\em planar} domains, through a formula on the fixed points of the involutions associated to the tilings. It was developed in 
\cite{Giraud,KTI,KTII}. We do not state it in detail here, since it does not play a role for the results in this letter. (The drum geometries describe general Gassmann-Sunada domains, not only those living in the real plane.)
\end{itemize}
\end{tabular}\right.$

One of the constructions is explained in \S\S \ref{mme}. (All the others are variations.)

Remarkably, the following result, the details of which can be found in \cite{KT-DRUMS}, tells us that no other examples can be constructed.

\begin{proposition}[\cite{KT-DRUMS}]
\label{propclass}
Length equivalent manifolds with FF, MAX, PAIR and INV which arise from the Gassmann-Sunada method,  {\em must} fall under one of the constructions of \cite{KT-DRUMS}, thus describing a precise algebraic form of these objects. 
\end{proposition}

Due to the nature of our constructions and properties, a connection with so-called {\em finite simple groups} occurs which seems, perhaps, rather surprising in this context. On the other hand, our properties define physically irreducible pairs of length equivalent manifolds | ``atoms'' of {\em general}
pairs of length equivalent manifolds, in that such a general pair of manifolds is patched up out of irreducible pairs | and that is precisely what simple groups are for general groups.  

Both in the constructions and in the proof of Proposition \ref{propclass}, a type of combinatorial geometry is present which is implicitly used throughout \cite{KT-DRUMS}. It shares a remarkable fixed points property with projective geometries, and is the subject of the next section.

\br{\rm
\begin{itemize}
\item[(a)]
In the known planar examples, the associated drum geometry, which is the main subject of this letter, allows an anti-isomorphism which switches $U$ and $V$. If one assumes the presence of such a mapping $\alpha$, one could put $W' = W^{\alpha}$ in MAX.
\item[(b)]
Note that if we would express PAIR through the existence of an outer automorphism of $G$ switching $U$ and $V$ (in which case the same number of tiles is trivially guaranteed), then it induces an anti-isomorphism as in (a), so (a) applies.  
\end{itemize}
}
\er

\medskip
\section{The geometry of drums}
\label{geomdrum}

A {\em drum geometry} or {\em D-geometry} is a pair $(\Gamma,A)$, in which $\Gamma$ is a rank $2$ geometry and $A$ a subgroup of the 
automorphism group $\Aut(\Gamma)$ (see Appendix \ref{geom} for formal definitions of these notions), such that the following properties are satisfied:

\begin{quote}
(TR)\quad $A$ acts transitively on both the points and the lines.
\end{quote}

\begin{quote}
(D)\quad An element $a \in A$ fixes at least one point of $\Gamma$ if and only if it fixes at least one line.
\end{quote}

If we use exponential notation, transitive on points, e.g., means that for any two points $x$ and $y$ there is an $a \in A$ such that
$x^a = y$.

A drum geometry $(\Gamma,A)$ is a {\em strong drum geometry} or {\em SD-geometry} if:
\begin{quote}
(SD)\quad An element $a \in A$ fixes precisely $m$ points of $\Gamma$ if and only if it fixes precisely $m$ lines.
\end{quote}

The main theorem of this section is the following.

\begin{theorem}[Geometry of drums]
\label{dgeom}
Any example of Gassmann-Sunada domains for the (strong) Dirichlet problem for the Laplacian arises  from a (strong) drum geometry.
Vice versa, every (strong) drum geometry gives rise to Gassmann-Sunada domains for the (strong) Dirichlet problem for the Laplacian. 
\end{theorem}

The connection works as follows. Let $T := (U,V,W)$ be a (strong) Gassmann-Sunada triple. Without loss of generality, we assume FF.
Define a geometry $\Gamma(T)$ as follows: points of $\Gamma(T)$ are 
elements of the left coset space $U/V = \{ uV \vert u \in U \}$; lines are elements of the left coset space $U/W = \{ uW \vert w \in W \}$. (Here, $uV := \{ u\cdot v \vert v \in V \}$.)
A point $aV$ is incident with the line $bW$ if and only if 
\begin{equation}
aV \cap bW \ne \emptyset. 
\end{equation}

%Now $U$ acts on this geometry by left translation, and one readily sees that $(U,\Gamma(T))$ is a drum geometry.

Conversely, let $(\Gamma,A)$ be a (strong) drum geometry, $x$  any point, and $L$ any line. Let $A_x$ be the stabilizer of $x$ in $A$, and $A_L$ the stabilizer of $L$ in $A$. Then the associated Gassmann-Sunada triple is $(A,A_x,A_L)$.\\

The details are below. In all these considerations, one only considers actions which are not sharply transitive, as these yield trivial results.\\

{\em Proof of Theorem \ref{dgeom}: details}.\quad
We use the notation of \S \ref{geomdrum}. We will prove the result for the strong version, because the general version follows trivially if the strong version is true.

First let $T := (U,V,W)$ be a strong Gassmann-Sunada triple with FF, and let $\Gamma(T) = (U/V,U/W,\bI)$ be the corresponding geometry, as in \S \ref{geomdrum}. Note that $U$ acts  on $\Gamma(T)$ by left translation, and obviously the action is transitive on both points and lines. For example, if $\widetilde{u}W$ is a line, and $u \in U$ arbitrary, then 
\begin{equation}
u(\widetilde{u}W) = (u\widetilde{u})W.
\end{equation}
So $U$ can be seen as an automorphism group of $\Gamma(T)$ (as it clearly preserves incidence), and one immediately sees that it is transitive on both points and lines.

We want to show that $(\Gamma(T),U)$ is a strong D-geometry.

The only thing we need is a simple fact about permutation group theory. Let $G$ be a group acting transitively on a set $X$  and suppose for convenience that no element of $G$ fixes every element of $X$ (that is, let $G$ be a transitive subgroup of the symmetric group $\texttt{Sym}(X)$ on $X$).
It is well known that this permutation group $(G,X)$ is equivalent to the permutation group $(G,G/G_x)$, where $x$ is any element in $X$, $G_x$ is, as before, the stabilizer of $x$ in $G$, and $G$ acts by left translation on the left cosets in $G/G_x$. For any element $g \in G$, denote by $\texttt{Fix}(g)$ the number of fixed points in $X$ (or $G/G_x$). Also, define $C_G(g) := \{ h \in G \vert gh = hg \}$, and put $g^G := \{ d^{-1}gd \vert d \in G \}$. Then
\begin{equation}
\label{fix}
\texttt{Fix}(g)\ = \ \frac{\vert C_G(g) \vert \cdot \vert g^G \cap G_x \vert}{\vert G_x \vert}.
\end{equation}

As $(U,V,W)$ is strong Gassmann-Sunada, we have by definition that 
\begin{equation}
\vert a^U \cap V \vert \ =\ \vert a^U \cap W \vert 
\end{equation}
for all $a \in U$. We know that $\vert V \vert = \vert W \vert$ by Theorem \ref{eqsizes}.  Now take any $u \in U$, and consider $u$ in its action of $(U,U/V)$ and $(U,U/W)$. Then inspecting $\texttt{Fix}(u)$ in both permutation representations, (\ref{fix}) gives us that this quantity is the same for both. 

In terms of the geometry $\Gamma(T)$, this means that the element $u$ fixes precisely the same number of points and lines | in other words, $(\Gamma(T),U)$ is a strong D-geometry.\\

Now let $(\Gamma,A)$ be a strong D-geometry, and let $(A,A_x,A_L)$ be as in \S \ref{geomdrum}. The fact that the geometry is strong means that 
for all $a \in A$, the number of fixed elements  in $(A,A/A_x)$ | that is, the number of fixed points | equals the number of fixed elements in $(A,A/A_L)$ | the number of fixed lines. (In particular, if $a$ is the identity element of $A$, it follows that $\vert A_x \vert = \vert A_L \vert$.) From (\ref{fix}) follows that 
\begin{equation}
\frac{\vert C_A(a) \vert \cdot \vert a^A \cap A_x \vert}{\vert A_x \vert}\ =\  \frac{\vert C_A(a) \vert \cdot \vert a^A \cap A_L \vert}{\vert A_L \vert}, 
\end{equation}
so that $\vert a^A \cap A_x \vert = \vert a^A \cap A_L \vert$ for all $a \in A$. 

Hence $(A,A_x,A_L)$ is strong Gassmann-Sunada. \eop \\

In terms of representations, we get:

\begin{corollary}[Representations and D-geometries]
Two faithful permutation representations $(G,G/U)$ and $(G,G/V)$ of the same degree (with $U$ and $V$ subgroups of the group $G$) have the same permutation character if and only if AC is satisfied for $(G,U,V)$ if and only if $\Gamma(G,U,V)$ is a strong drum geometry.
\eop
\end{corollary}

\subsection{Very important example (well known)}
\label{VE1}

Let $k$ be any field (one may take $k = \mathbb{C}$, for instance), let for any positive integer $n \geq 3$,
$\PGL_n(k)$ be the corresponding projective general linear group. Its elements are nonsingular $(n \times n)$-matrices with entries in $k$, and 
they are defined up to scalars. (The group operation is just matrix multiplication.) It is well known that $\PGL_n(k)$ has a natural action on 
the $(n - 1)$-projective space $\PG(n - 1,k)$: after having chosen homogeneous coordinates in that space, the image of a point $\mathbf{x}$ (written in coordinates) under an element 
$A \in \PGL_n(k)$ is 
\begin{equation}
A \cdot \mathbf{x}^{\text{T}}.
\end{equation} 

In a natural way all elements of $\PGL_n(k)$ also act on subspaces of $\PG(n - 1,k)$.

In this section, if $\kappa$ is a subspace of $\PG(n - 1,k)$, then $\PGL_n(k)_{\kappa}$ is the stabilizer in $\PGL_n(k)$ of $\kappa$. 

The following facts are well known and easy to prove.

\begin{itemize}
\item[{\bf Fact I}]
$\PGL_n(k)$ acts transitively on the points and hyperplanes  (= $(n - 2)$-dimensional subspaces) of $\PG(n - 1,k)$.
\end{itemize}

\begin{itemize}
\item[{\bf Fact II}]
An element of $\PGL_n(k)$ fixes at least one point of $\PG(n - 1,k)$ if and only if it fixes at least one hyperplane.
\end{itemize}

Hence, if we define the geometry $\Gamma$ to have as points the points of $\PG(n - 1,k)$ and as lines the $(n - 2)$-subspaces, incidence being 
inherited from the projective space, we have that $(\Gamma,\PGL_n(k))$ is a D-geometry.

\subsection{Very important example (well known)}
\label{VE2}

Now let $k$ be a {\em finite} field, so that the number of subspaces of $\PG(n - 1,k)$ is finite, and the size of $\PGL_n(k)$ is as well.  
The following fact is not hard to prove (see also \S\S \ref{3.4}), but important: 

\begin{itemize}
\item[{\bf Fact III}]
An element of $\PGL_n(k)$ fixes $c$ points of $\PGL(n - 1,k)$ if and only if it fixes $c$ hyperplanes.
\end{itemize}

So with the same notation as above, we have that $(\Gamma,\PGL_n(k))$ is a {\em strong} D-geometry. Up to natural equivalence, 
we get the corresponding Gassmann-Sunada triple as before: let $x$ be any point and $Y$ be any hyperplane; then 
\begin{equation}
(\PGL_n(k),{\PGL_n(k)}_x,{\PGL_n(k)}_Y)
\end{equation}
is the required triple.\\

All these examples have FF, MAX, PAIR. All the examples of this type that give rise to {\em planar} domains were determined in \cite{KTI,KTII,KTIII}, and correspond to the presently known planar examples. They are precisely the examples with INV.

\subsection{Many more examples (\cite{KT-DRUMS})}
\label{mme}

We will define the class $\mD$ of ``Type I triples'' of \cite{KT-DRUMS}.

Let $(\Gamma,A)$ be any D-geometry, and note again that strong D-geometries are also D-geometries. One may take any of the previous examples, for instance. 

Let $(A,A_x,A_L)$ be ``the'' corresponding Gassmann-Sunada triple (which has FF). 

New Gassmann-Sunada triples are constructed through the following data: 
\begin{equation}
\mathcal{E} = (A,A_x,A_L,n,T)
\end{equation}
where 
\begin{itemize}
%\item
%$L$ and $L'$ are not conjugate (and in particular $L \ne L'$),
%\item
%$(S,L,L')$ is an EC-triple with FF (we call $(S,L,L')$ the {\em base} of $\mE$), 
\item
$n \in \mathbb{N}^\times$ is any positive integer, and 
\item
$T$ is any transitive subgroup of the symmetric group $\texttt{S}_n$ acting naturally on $\{1,\ldots,n\}$. 
\end{itemize}

Define a triple $(A \wr T,A_x \wr T,A_L \wr T)$, where ``$\wr$'' denotes the wreath product (see \cite{KT-DRUMS} for details). Note that 
$\vert A  \wr T\vert = \vert A\vert^n\vert T\vert$,   $\vert A_x \wr T\vert = \vert A_x\vert^n\vert T\vert$ and  $\vert A_L \wr T\vert = \vert A_L\vert^n\vert T\vert$.

The following fundamental properties are proven in \cite{KT-DRUMS}.
\begin{itemize}
\item
$(A \wr T,A_x \wr T,A_L \wr T)$ is a Gassmann-Sunada triple  which satisfies FF.
\item
If $A_x$ and $A_L$ are not conjugate, then $A_x \wr T$ and $A_L \wr T$ also are not conjugate.
\item
If $(A_,A_x,A_L)$ has MAX, $(A \wr T,A_x \wr T,A_L \wr T)$ also has MAX. 
\item
If $(A,A_x,A_L)$ has PAIR, $(A \wr T,A_x \wr T,A_L \wr T)$ also has PAIR.
\end{itemize}

In other words, the newly constructed domains are also length equivalent, nontrivial, and irreducible!

As each finite group is a transitive group, the generality of the construction is clear. The paper \cite{KT-DRUMS} contains two variations of this construction (``Type II triples'' and ``Type III triples''), and one of the main results in \cite{KT-DRUMS} reads that any Gassmann-Sunada triple  that satisfies FF, MAX and PAIR 
is of one of these types. So \cite{KT-DRUMS} contains a complete classification of physically irreducible Gassmann-Sunada drum pairs.

We end this subsection with the following result, taken from \cite{KT-DRUMS}.

\begin{theorem}[Unbounded index]
For each $N \in \mathbb{N}$, there exist Gassmann-Sunada triples $(G,U,V)$ with MAX and FF, such that
\begin{equation}
\frac{\vert G\vert}{\vert U\vert} = \frac{\vert G\vert}{\vert V\vert} > N.
\end{equation}
\end{theorem}

\subsection{Still more examples:  (transitive) symmetric designs}
\label{3.4}

As we now know, any strong Gassmann-Sanada pair of domains is associated to a strong D-geometry, and the observation reverses. To show the power of this link, we will define a very large class of geometries which have SD. 

A {\em symmetric $2$-$(v,k,\lambda)$ design} consists of two kinds of objects usually called {\em points} and {\em blocks} (instead of lines), such that any two different points are contained in $\lambda$ blocks and any two different blocks meet in $\lambda$ points; any block has $k$ points and any point is in $k$ blocks; there are $v$ points and $v$ blocks. This class of combinatorial geometries is huge, as its literature shows (see e.g. \cite{symmdes}). 
Our interest is based on the following standard result.

\begin{theorem}[\cite{AB}]
\label{symmfix}
Let $\Gamma$ be a symmetric $2$-$(v,k,\lambda)$ design, and let $\alpha$ be an automorphism of $\Gamma$. Then $\alpha$ fixes precisely $c \in \mathbb{N}$ points if and only if $\alpha$ fixes precisely $c$ blocks. In other words, (SD) is satisfied w.r.t. any automorphism. Moreover, any automorphism group is transitive on points if and only if it is on blocks.
\end{theorem}

\begin{corollary}
\label{cor3.4}
Let $\Gamma$ be a symmetric $2$-$(v,k,\lambda)$ design and let $A$ be a point-transitive automorphism group. Then $(\Gamma,A)$ is a strong D-geometry.
\eop
\end{corollary}

All the examples in \S\S \ref{VE1} and \S\S \ref{VE2} fall under this general construction, but they are just the tip of the iceberg.

In the forthcoming paper \cite{desdrum}, we will analyze the examples of this subsection in further detail. We give one interesting example.

Let $V = V(2n,\F_2)$ a $2n$-dimensional vector space over the field with two elements $\F_2$. Let $q: V \mapsto \F_2$ be a nonsingular quadratic form in $2n$ variables over $\F_2$; then $D := \{ v \ \vert\ q(v) = 1 \}$ is a so-called {\em difference set} in the additive group $(V, +)$, with parameters $(2^{2n},2^{2n - 1} + \epsilon 2^{n - 1},2^{2n - 2} + \epsilon 2^{n - 1})$. The sign $\epsilon \in \{-1,+1\}$ depends on the choice of the quadratic form.  One constructs a symmetric $2$-$(2^{2n},2^{2n - 1} + \epsilon 2^{n - 1},2^{2n - 2} + \epsilon 2^{n - 1})$ design $\mS^{\epsilon}(2n)$ from this data by letting the points be the elements of $V$, and the lines be the sets $\{ a + D\ \vert\ a \in V \}$ (incidence is symmetrized containment). The automorphism group of this design acts $2$-transitively on the points and blocks (cf. \cite{Kantor} for more details).

\bt
Let $A$ be the full automorphism group of $\mS^{\epsilon}(2n)$, let $x$ be a point and $L$ be a line.
Then $(A,A_x,A_L)$ is a strong Gassmann-Sunada triple satisfying FF, MAX and PAIR. In general, INV is not satisfied. 
\et

{\em Proof}.\quad
By Corollary \ref{cor3.4}, the drum geometry $(\mS^{\epsilon}(2n),A)$ is strong, so that $(A,A_x,A_L)$ is strong Gassmann-Sunada. Hence the triple satisfies PAIR. By definition of $A$, FF is also true. Since $A$ acts $2$-transitively on the left coset space $A/A_x$, it also act primitively on it, so $A_x$ is maximal in $A$. The same holds for $A_L$, whence MAX is satisfied. 

Finally, all the strong Gassmann-Sunada triples $(U,V,W)$ where $U$ acts $2$-transitively on $U/V$ and $U/W$, and which satisfy INV, are classified in \cite{JSKT}.
\eop

\subsection{Dihedral groups}

The following example is certainly well known (and trivial), but we introduce it  to illustrate the ease of passing to D-geometries. Consider a dihedral group $A := \mathbf{D}_{n}$ of order  $2n$, $n \in \mathbb{N} \setminus \{ 0, 1\}$; it is a group with representation $\Big\langle u, v\ \vert\ u^n = v^2 = (uv)^2 = 1 \Big\rangle$. We see it as the automorphism group of the point-line geometry $\Gamma$ defined by the vertices and sides of a regular $n$-gon embedded in $\mathbb{R}^2$. Then any element $a \in A^\times$ has 
either $0$ points and $0$ lines; or $1$ point and $1$ line (which can only happen if $n$ is odd); or $2$ points and $0$ lines (which can only happen if $n$ is even); or $0$ points and $2$ lines (still in the $n$ even case). So if $n$ is odd, with $x$ any point of $\Gamma$ and $Y$ any line, $(A,A_x,A_Y)$ is a strong Gassmann-Sunada triple (with $\vert A_x \vert = \vert A_Y \vert = 2$). However, $A_x$ and $A_Y$ are obviously conjugate for each choice of $x$ and $Y$.

\subsection{Examples coming from reductive groups}

In this subsection, I describe an abstract construction of strong Gassmann-Sunada triples which might well be the most general one described 
in the literature. It was communicated to me by Dipendra Prasad, and is mentioned in his paper \cite{Prasad}. We will not give precise definitions
for all the notions used, but suggest references instead. 

Let $G$ be a reductive group defined over a finite field $k = \mathbb{F}_q$. Let $P_1$ and $P_2$ be parabolic subgroups such that the Levi subgroups $\mL_1$ and $\mL_2$ of respectively $P_1$ and $P_2$ are conjugate in $G$. Then $(G(k),P_1(k),P_2(k))$ is a strong Gassmann-Sunada triple \cite{Prasad}. This already comes with a geometric interpretation. Each reductive group has {\em BN-pairs} (as defined in \cite{Tits}), and up to conjugation, one can assume that there is a BN-pair $(B,N)$ in $G(k)$ of which $P_1$ and $P_2$ are parabolics. (In the context of
almost conjugacy, we may conjugate the subgroups $U$ and $V$ in a (strong) Gassmann-Sunada triple $(U,V,W)$; cf. equation (\ref{eqisom}) in \S \ref{isom}.) So there is a Tits building $\mB$ for which left cosets $gP_1$ and  $hP_2$ represent elements of the types corresponding to the parabolics.  Passing to the associated D-geometry $(\Gamma,G(k))$, what we do is: 
\begin{itemize}
\item[$\star$]
precisely single out the elements of these types, 
\item[$\star$]
call the elements ``points'' and ``lines,''
\item[$\star$]
endow the incidence of the building on this new geometry. 
\end{itemize}

The property about Levi subgroups implies property (SD) for this D-geometry. Specializing to the case $G = \PSL_{n + 1}(q)$ acting naturally on projective space $\PG(n,q)$, and letting $P_1$ and $P_2$ correspond to points and hyperplanes, we obtain the construction of section \ref{VE2}.

\medskip
\section{Strong D-geometries and transplantation (new proof of an old theme)}

For all strong D-geometries, there is an easy way to construct a transplantation matrix between the two respective linear representations, without effectively passing to linear representations (cf the introduction). We will handle this aspect in a general context using a standard matrix representation method for point-line geometries.

Let $\Gamma = (\mP,\mL,\I)$ be any point-line incidence geometry. We suppose it is finite, meaning that the product $\vert \mP \vert \times \vert \mL \vert$ is finite. Let $\mu = \vert \mP \vert$ and $\nu = \vert \mL \vert$. Choose and fix a labeling $x_1,\ldots,x_\mu$ and $Y_1,\ldots,Y_\nu$ of both the point and line set. Now we define a $(\mu \times \nu)$-matrix as follows:
\begin{equation}
\begin{cases}
a_{ij} = 1 &\mbox{if}\ x_i \I Y_j,\\
a_{ij} = 0 &\mbox{otherwise}. 
\end{cases}
\end{equation}

Now let $\alpha$ be an arbitrary automorphism of $\Gamma$; then $\alpha$ induces a permutation of $\mP$ and a permutation of $\mL$. Define the $(\mu \times \mu)$-matrix $P$ as follows:
\begin{equation}
\begin{cases}
p_{ij} = 1 &\mbox{if}\ \alpha(x_i) = x_j,\\
p_{ij} = 0 &\mbox{otherwise}. 
\end{cases}
\end{equation}

Similarly we define a $(\nu \times \nu)$-matrix $L$ with respect to the line set. Note that $P$ and $L$ are permutation matrices.

\bp
\label{propcomm}
We have that $PA = AL$.
\ep
{\em Proof}.\quad
Consider the element $(PA)_{ij} = \sum_{\rho = 1}^\mu p_{i\rho}a_{\rho j}$; then this entry equals $a_{u j}$, where $p_{iu} = 1$, that is,
 $\alpha(x_i) = x_u$. Similarly, we have $(AL)_{ij} = \sum_{\delta = 1}^\nu a_{i\delta}l_{\delta j}$; this entry equals $a_{i v}$, where $l_{vj} = 1$, meaning, $\alpha(Y_v) = Y_j$. Now $(PA)_{ij} = 1$ if and only if $a_{uj} = 1$ (with $\alpha(x_i) = x_u$) if and only if $x_u \I Y_j$ if and only if $\alpha^{-1}(x_u) \I 
 \alpha^{-1}(Y_j)$ if and only if $x_i \I Y_v$ if and only if $a_{iv} = 1$. \eop \\

Note that $A$ is defined over $\mathbb{Z}$ (as it only contains entries equal to $0$ or $1$). \\

Since no additional assumptions were made on the geometry $\Gamma$, the identity in Proposition \ref{propcomm} does not express transplantability (yet), although it looks like it. Now suppose that $\mu = \nu$, and that
$A \in \GL_{\mu}(k)$ for some field $k$ of characteristic $0$, or for $k = \mathbb{Z}$; then $A^{-1}PA = Q$, and so $\mathrm{trace}(P) = 
\mathrm{trace}(Q)$ for all $\alpha \in K$. This means precisely that each $\alpha \in K$ fixes the same number of points and lines. We have proved the following proposition (stated slightly more generally):

\bp
\label{propgeomtran}
Let $\Gamma = (\mP,\mL)$ be a finite point-line geometry, where $\vert \mP \vert = \vert \mL \vert$, and let $A$ be an incidence matrix for $\Gamma$.
If $A$ is invertible over $\mathbb{Z}$, or some field of characteristic $0$, then for any element $\alpha$ in $\Aut(\Gamma)$, we have that $PA = AQ$ (using the 
notation of above). In particular, if $K \leq \Aut(\Gamma)$ is any (point- or line-) transitive automorphism group, and $(x,Y)$ is any point-line pair 
of $\Gamma$, then $(K,K_x,K_Y)$ is strong Gassmann-Sunada.  
\ep

Note that the choice of labeling of the points and lines does not influence the result; such changes amount to passing from $A$ to an equivalent matrix. 

Proposition \ref{propgeomtran} and the discussion before that, motivate the next definition.
\begin{quote}
$\Big[${\bf Super strong}$\Big]$.\quad
A point-line geometry $\Gamma = (\mP,\mL)$ with $\vert \mP \vert = \vert \mL \vert$ is {\em super strong} if for all $\alpha \in \Aut(\Gamma)$, we have that $\alpha$ fixes equally as many points as lines. This is equivalent to saying that any (= some) incidence matrix be non-singular over $\mathbb{R}$. 
\end{quote}

\br{\rm
\begin{itemize}
\item[{\rm (a)}]
If $(\Gamma,K)$ is a strong D-geometry (with associated Gassmann-Sunada triple $(K,K_x,K_Y)$), $\mu = \nu$, and for every element $\alpha \in K$, $\alpha$ induces an element in $(K,K/K_x)$ and $(K,K/K_Y)$. Representing the elements as permutation matrices $P$ and $L$ in the usual way, we obtain
\begin{equation}
PA = AL,
\end{equation}
where for each $\alpha$, the traces of $P$ and $Q$ are the same. This does not give any guarantee that $A$ is invertible over $\mathbb{R}$. It would be interesting to have examples at hand of strong D-geometries which are not super strong, i. e., for which the incidence matrices are singular.
\item[{\rm (b)}]
Still working in the same situation, note that in any case there is an element $\widehat{A} \in \GL_{\mu}(\C)$ such that $P\widehat{A} = \widehat{A}Q$ for all $\alpha \in K$. Now consider the set $\mA$ of all such matrices $\widehat{A}$: if $\mA$ contains an element $\widetilde{A}$ with only $0$ or $1$ as entries, and with only zeroes on the  diagonal, then one could define a point-line geometry $\Gamma(\widetilde{A})$ with $\widetilde{A}$ as incidence matrix. It is easy to show that each $\alpha \in K$ naturally induces an automorphism of $\Gamma(\widetilde{A})$; since $\widetilde{A}$ is not singular, the geometry is super strong.
\item[{\rm (c)}]
Almost no cases are known in which $\widehat{A}$ is an element of $\GL_{\mu}(\mathbb{Z})$, i.e., for which $\vert \mathrm{det}(\widehat{A}) \vert = 1$; see 
the recent paper of Prasad \cite{Prasad}. Such Gassmann-Sunada triples give rise to nonisomorphic projective curves  over fields $k$ with isomorphic Jacobians (so nonisomorphic curves with the same motive) \cite{Prasad}.   I do not know what the (incidence) geometrical meaning is of this property. In any case, it is well known that symmetric $2$-designs have incidence matrices whose determinant is different from $1$. 
\end{itemize}
}\er

%\subsection*{Checking the Gassmann-Sunada property locally}

%\subsection*{Corollary: checking whether a D-geometry is super strong}

\medskip
\section{FF and D-geometries}

The properties FF, MAX, PAIR and INV are easily interpreted on D-geometries. In this section, we have a look at FF.

Let $(U,V,W)$ be a Gassmann-Sunada triple which does not satisfy FF. So there is a normal subgroup $N$ of $U$ which is also contained in $V$ and $W$. Let $n \in N$. Then for all $u \in U$, we have that $n(uV) = un'V = uV$ for some $n' \in N$, and similarly we have that $n(uW) = uW$. In other 
words, $N$ acts trivially on the points and lines of the associated D-geometry. 

Vice versa, let $(\Gamma,U)$ be a D-geometry, and suppose $U$ does not act faithfully on the points or lines of $\Gamma$. Let $N$ be the subgroup of $U$ which fixes all lines of $\Gamma$. As we identify a line with the set of points incident with it, this is also the subgroup of $N$
of which each element fixes each line of $\Gamma$. Then $(U,V,W)$ does not satisfy FF, but $(U/N,V/N,W/N)$ does.

In the rest of this paper, we only consider D-geometries with FF.

\medskip
\section{Constructions via subgroups and overgroups; Isomorphism testing}
\label{isom}

It might be desirable to introduce a notion of {\em isomorphism} for Gassmann-Sunada triples, and also of D-geometries, in order to 
see (eventually) when, for instance, two given  Gassmann-Sunada triples give isomorphic topological or arithmetic objects. This is in fact quite tricky:
suppose for example that $(U,V,W)$ is a (strong or not) Gassmann-Sunada triple with FF, and construct the associated D-geometry $(\Gamma,U)$.  Then $V$ is the stabilizer in $U$ of some point $x$, and $W$ is the stabilizer of some line $Y$. Now let $x'$ be another point, and let $Y'$ be another line. Then $(U,U_{x'},U_{Y'})$ is also (strong or not) Gassmann-Sunada; this is easily seen as follows. First note that there are elements $u, v \in U$ such that $V^u = U_{x'}$ and $W^v = U_{Y'}$ (here, we use exponential notation, so that $V^u = u^{-1}Vu$, e.g.).  
Let $\mG$ be any conjugacy class in $U$; then
\begin{equation}
\label{eqisom}
\begin{cases}
\vert \mG \cap V \vert = \vert \mG \cap V^u \vert = \vert \mG \cap U_{x'} \vert \\
\vert \mG \cap W \vert = \vert \mG \cap W^v \vert = \vert \mG \cap U_{Y'} \vert.
\end{cases}
\end{equation}

The $U$-actions $(U,U/U_x)$ and $(U/U_{x'})$ are isomorphic, and so are 
$(U,U/U_Y)$ and $(U,U/U_{Y'})$. So the triples $(U,V,W)$ and $(U,U_{x'},U_{Y'})$ should be considered to be isomorphic.

Based on these observations, 
we say that two
Gassmann-Sunada triples $(U,V,W)$ and $(U',V',W')$ are {\em isomorphic} if the exists a group isomorphism $\phi: U \mapsto U'$, and 
elements $u, v \in U$, such that
${(V^\phi)}^u = V'$ and ${(W^\phi)}^v = W'$. On the level of geometries, this notion translates as follows: 

\begin{quote}
D-geometries $(\Gamma,U)$ and 
$(\Gamma',U')$ are isomorphic if there exists a point line pair $(x,Y)$ in $\Gamma$ and $(x',Y')$ in $\Gamma'$ such that the 
Gassmann-Sunada triples $(U,U_x,U_Y)$ and $(U',U'_{x'},U'_{Y'})$ are isomorphic. 
\end{quote}

This is the same as asking that they are isomorphic for 
any choice of such pairs.

\medskip
\section{Constructing new (strong) Gassmann-Sunada triples from olds ones}

Without further mention, we suppose that each Gassmann-Sunada triple in this section satisfies FF.

Suppose $(U,V,W)$ is a (strong) Gassmann-Sunada triple, and construct the geometry $(\Gamma,U)$. Then we know that (D) or (SD) is satisfied, depending on whether $(\Gamma,A)$ is strong or not. Now let $X$ be a subgroup of $U$ which is still transitive on points and lines. Then (D) or (SD) is still satisfied for $(\Gamma,X)$, and we obtain a (strong) Gassmann-Sunada triple $(X,X_y,X_L)$, where $x$ is any point and $L$ is any line. 

\bp
Suppose $(U,V,W)$ is a (strong) Gassmann-Sunada triple, and construct the geometry $(\Gamma,U)$. Let $X$ be a subgroup of $U$ which is still transitive on points and lines. Let $y$ be any point and $L$ be any line.
We have that $(X,X_y,X_L)$ is (strong) Gassmann-Sunada.
\ep

What about overgroups? It is easy to see that if $(U,V,W)$ is a (strong) Gassmann-Sunada triple, and $\overline{U} \geq U$ is any group containing $U$, then 
$(\overline{U},V,W)$ is also (strong) Gassmann-Sunada. Such examples are usually not considered as being ``new.'' There is a much more subtle 
construction; it is contained in the next theorem.

%If $\Gamma$ is a symmetric $2$-design, part (b) of the next theorem gives a complete answer.

\bt
\label{enl}
\begin{itemize}
\item[{\rm (a)}]
Suppose that $\Gamma$ is super strong, and that $A \leq \Aut(\Gamma)$ is transitive on points and lines. If $A \leq B \ne A \leq \Aut(\Gamma)$, and 
$(x,Y)$ is any point-line pair, then $(B,B_x,B_Y)$ is strong Gassmann-Sunada. 
\item[{\rm (b)}]
Suppose $(U,V,W)$ is a strong Gassmann-Sunada triple, and construct the geometry $(\Gamma,U)$;
suppose furthermore that $\Gamma$ is a symmetric $2$-design. Let $\widetilde{U}$ be a subgroup of $\Aut(\Gamma)$ in which $U$ is contained. Then $(\widetilde{U},\widetilde{U}_y,\widetilde{U}_L)$ is strong Gassmann-Sunada.
\end{itemize}
\et
{\em Proof}.\quad
Since $\Gamma$ is super strong and $A$ is transitive on points and lines, part (a) easily follows. 
In part (b), by Theorem \ref{symmfix}, (SD) is satisfied for $U$, whence the result. \eop \\

So the structure of the geometry yields extra information to construct ``new'' Gassmann-Sunada triples. But what does ``new'' mean? Can triples arising as in (a) or (b) yield different isospectral manifolds than the ones obtained by using $(U,V,W)$? Since there are many different methods to construct such manifolds from a given triple, the manifold aspect is not clear to us, and it would be desirable to have a precise connection between isomorphisms on the level of D-geometries, and isomorphisms (= isometries) on the level of manifolds.

A geometry could be D/SD for several nonisomorphic groups at the same time. For instance, use the notation of Theorem \ref{enl}(a); then on the geometrical level, we have:
\begin{equation}
(B,B/B_x,B/B_L)\ \cong\ (A,A/A_x,A/A_L)\ \cong \Gamma.
\end{equation}

On the other hand, $B$ might give extra information; let $A$ be, e.g., a transitive but not $2$-transitive group (on points and lines), and let $B$ be $2$-transitive (on points and/or lines). Other example: consider a strong Gassmann-Sunada triple $(U,V,W)$ (with FF), pass to the drum geometry $(\Gamma,U)$ which we suppose to be super strong, and let $U \leq \widetilde{U} \leq \Aut(\Gamma)$ with $U \ne \widetilde{U}$. Then $(\widetilde{U},\widetilde{U}_x,\widetilde{U}_L)$ is strong Gassmann-Sunada, and FF by construction. Since $\widetilde{U} \leq \Aut(\Gamma)$, points and lines are mapped to points and lines, respectively, and so $\widetilde{U}_x$ and $\widetilde{U}_L$ stay non-conjugate. If we would allow dualities in an extended
notion of automorphism group, they still could become conjugate. We explain this in some more detail below.

First, define a {\em duality} $\rho$ of a point-line geometry $\Delta = (\mP,\mL,\I)$ to be a permutation of $\mP \cup \mB$ such that $\rho(\mP) = \mL$, $\rho(\mL) = \mP$, and $x \I Y$ if and only if $\rho(x) \I \rho(Y)$. Let $\Gamma$ be a strong D-geometry with respect to the automorphism 
group $A$, so that we obtain a strong Gassmann-Sunada triple $(A,A_x,A_L)$; then by the discussion prior to the statement of Theorem \ref{enl}, we have that $(\Aut(\Gamma),A_x,A_L)$ is also strong Gassmann-Sunada. If $\delta$ is a duality of $\Gamma$, it follows that $(\langle \Aut(\Gamma),\delta \rangle,A_x,A_L)$ is also strong Gassmann-Sunada. On the other hand, $V^\delta = W$ and $W^\delta = V$, so $V$ and $W$ are conjugate in $\langle \Aut(\Gamma),\delta \rangle$. 

Now let $(U,V,W)$ be strong Gassmann-Sunada, and suppose  that $V^u = W$ for some element $u \in U$. Suppose moreover that $u^2$ normalizes $V$, so that $W^u = V$. (Special example: the case $u^2 = \id$.) Then obviously $u$ induces a duality of the geometry $\Gamma(U,V,W)$. \\

%\bigskip
\newpage
\appendix

\section{The Gordon--Webb--Wolpert example}
\label{count}

The most famous examples of isospectral but not isometric planar domains are  shown in Figure \ref{celebrated}, and are essentially taken from \cite{GWW}.
\begin{Figure}
\begin{center}
\includegraphics[width=6cm]{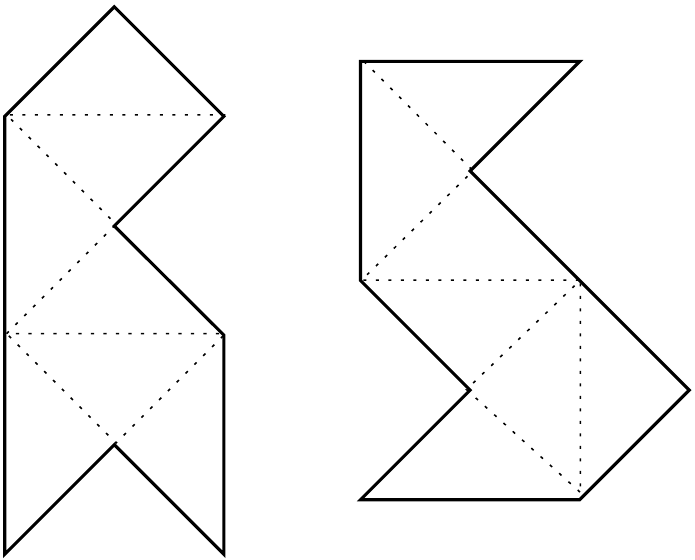}
\end{center}
\label{celebrated}
\end{Figure}

Both drums are constructed by unfolding the same triangle seven times, along the dotted lines (according to a transplantation rule).

The $7$ tiles on the left correspond to the points (or lines) of $\PG(2,2)$, the tiles on the right to lines (or points) of $\PG(2,2)$. So in this case, the 
strong D-geometry is $\PG(2,2)$ itself (or more precisely, $(\PG(2,2),\PGL_3(2))$).

The corresponding Gassmann-Sunada triple has FF, MAX, PAIR and INV.

\section{The Gassmann-Sunada method}
\label{GS}

In this appendix, we describe the Gassmann-Sunada method, and refer to \cite{Gassman,Su} for original work.

Let $H, H'$ be subgroups of a finite group $G$. Then $H$ and $H'$ are {\em almost conjugate} in $G$ if for all $g \in G$, we have that $\vert g^G \cap H \vert = \vert g^G \cap H' \vert$. Here, $g^G = \{ h^{-1}gh \vert h \in G \}$.

The following properties are easily shown to be equivalent. \begin{itemize}
\item[(a)]
$H$ and $H'$ are almost conjugate (AC) (that is, $(G,H',H')$ is by definition a strong {\em Gassmann-Sunada triple}).
\item[(b)]
There exists a bijection $B: H \longrightarrow H'$ such that for each $h \in H$, $B(h)$ is conjugate to $h$ in $G$.
%\item[(c)]
\end{itemize}

The next theorem is easy to prove, and well known.

\begin{theorem}
\label{eqsizes}
If $H$ and $H'$ are AC in $G$, then $\vert H\vert = \vert H'\vert$.
\end{theorem}

\begin{theorem}[AC and Isospectrality, \cite{Su}]
\label{ACIso}
Let $\M$ be a compact Riemannian manifold, $G$ a group, and $H$ and $K$ almost
conjugate subgroups of $G$. Then if $\pi_1(\M)$ admits a homomorphism onto $G$, the covers $\M_H$
and $\M_K$ associated to the pullback subgroups of $H$ and $K$ are isospectral.
\end{theorem}

We say that subgroups $H$ and $H'$ of a finite group $G$ are {\em elementwise conjugate} if for each $g \in G$ we have that
\begin{equation}
g^G \cap H \ne \emptyset \ \   \Longleftrightarrow\ \      g^G \cap H' \ne \emptyset
\end{equation}

The following properties are  equivalent.
\begin{itemize}
\item[(i)]
$H$ and $H'$ are almost elementwise conjugate (EC).
\item[(ii)]
Each element of $H'$ is conjugate to an element of $H$ and each element of $H$ is conjugate to some element of $H'$.
\end{itemize}

For $\M$ a compact Riemannian manifold, let $\mathcal{L}(\M)$ be the set of lengths of closed geodesics of $\M$ (that is, $\mathcal{L}(\M)$ is the {\em length spectrum} without multiplicities). Then two such manifolds $\M$ and $\widetilde{\M}$ are {\em length equivalent} if $\mathcal{L}(\M) = \mathcal{L}(\widetilde{\M})$.

Expressed in terms of isospectrality-related properties, we have the following version of Theorem \ref{ACIso}.

\begin{theorem}[EC and Length Equivalence \cite{Lein}]
\label{ECLE}
Let $\M$ be a compact Riemannian manifold, $G$ a group, and $H$ and $K$ elementwise
conjugate subgroups of $G$. Then if $\pi_1(\M)$ admits a homomorphism onto $G$, the covers $\M_H$
and $\M_K$ associated to the pullback subgroups of $H$ and $K$ are length equivalent.
\end{theorem}

The following is obvious.
\begin{observation}
\label{EXAC}
Let $(G,H,H')$ be a strong Gassmann-Sunada triple. Then $H$ and $H'$ are also elementwise conjugate in $G$.
\end{observation}

Triples $(G,H,H')$ in which $H$ and $H'$ are EC, are called (``ordinary'') Gassmann-Sunada triples.

So in general we have
%\begin{equation}
%\mathrm{AC} \ \ {\huge \mapsto}\ \ \mathrm{EC}.
%\end{equation}

\begin{center}
\begin{tikzpicture}[>=angle 90,scale=2.2,text height=1.5ex, text depth=0.25ex]
%%First place the nodes
\node (a0) at (0,3) {AC};
\node (a1) [right=of a0] {EC};
\draw[->,font=\scriptsize,thick]
(a0) edge node[auto] {implies} (a1);
\end{tikzpicture}
\end{center}

The converse is not necessarily true | one cannot even conclude that $H$ and $H'$ have the same order (cf. the PAIR property)!
For example, let $\texttt{A}_n$ be the alternating group on $n$ letters ($n \geq 3$), and consider commuting involutions (transpositions) $\alpha$ and $\beta$.
Then $H = \{ \id,\alpha \}$ and $H' = \{ \id,\alpha,\beta,\alpha\beta \}$ are EC (and certainly not AC, as $\vert H \vert \ne \vert H'\vert$).

\section{Notions of geometry}
\label{geom}

A {\em rank $2$ geometry} is a triple $\Gamma = (\mP,\mL,\bI)$, in which $\mP$ and $\mL$ are disjoint sets called, respectively, ``points'' and ``lines,''
and where $\bI \subseteq (\mP \times \mL) \cup (\mL \times \mP)$ is a binary irreflexive symmetric relation on $\mP \cup \mL$, which expresses when a point is 
incident with a line. 

An {\em automorphism} $\alpha$ of $\Gamma$ is a permutation of $\mP \cup \mL$ which preserves $\mP$ and $\mL$, with the additional property that $x \bI Y$ if and 
only if $x^{\alpha} \bI Y^{\alpha}$ | in other words, it preserves incidence. The set of all automorphisms of the geometry $\Gamma$ forms a group under 
compostion, which is denoted by $\Aut(\Gamma)$. 

A standard reference on incidence geometry is \cite{HB}.

%\section{BN-Pairs}
%\label{BN}

%\bigskip

%\end{multicols}

\end{document}